\documentclass[11pt]{article}

\usepackage{amsfonts,epsf,amsmath,amssymb}
\usepackage{epsfig}
\usepackage{latexsym,multirow,rotating}
\usepackage{pgf,tikz,subfigure}
\usetikzlibrary{decorations.markings}

\newtheorem{theorem}{\bf Theorem}[section]

\newtheorem{lemma}[theorem]{\bf Lemma}
\newtheorem{proposition}[theorem]{\bf Proposition}

\newtheorem{definition}[theorem]{\bf Definition}

\newtheorem{fact}[theorem]{\bf Fact}

\newcommand{\proof}{\noindent{\bf Proof.\ }}
\newcommand{\qed}{\hfill $\square$ \bigskip}

\begin{document}

\title{ANTIPARALLEL $d$-STABLE TRACES AND A STRONGER VERSION OF ORE PROBLEM}

\author{
Jernej Rus \\
Faculty of Mathematics and Physics, University of Ljubljana \\
Jadranska 19, 1000 Ljubljana, Slovenia \\
jernej.rus@gmail.com
}
\date{\today}
\maketitle

\begin{abstract}

In $2013$ a novel self-assembly strategy for polypeptide nanostructure design which could lead to significant developments in biotechnology was presented in [Design of a single-chain polypeptide tetrahedron assembled from coiled-coil segments, Nature Chem. Bio. 9 (2013) 362--366]. It was since observed that a polyhedron $P$ can be realized by interlocking pairs of polypeptide chains if its corresponding graph $G(P)$ admits a strong trace. It was since also demonstrated that a similar strategy can also be expanded to self-assembly of designed DNA [Design principles for rapid folding of knotted DNA nanostructures, Nature communications 7 (2016) 1--8.]. In this direction, in the present paper we characterize graphs which admit closed walk which traverses every edge exactly once in each direction and for every vertex $v$, there is no subset $N$ of its neighbors, with $1 \leq |N| \leq d$, such that every time the walk enters $v$ from $N$, it also exits to a vertex in $N$. This extends C. Thomassen's characterization [Bidirectional retracting-free double tracings and upper embeddability of graphs, J. Combin. Theory Ser. B 50 (1990) 198--207] for the case $d = 1$.
\end{abstract}


\newpage
%%%%%%%%

%%%%%%%%%%%%%%%%%
\section{Introduction}

Nucleic acids and proteins are programmable polymers where the sequence of their constituent monomers in a chain defines the self-assembly of proteins and nucleic acids into defined tertiary structures that underlies their function. Natural globular proteins are stabilized by a complex balance of a large number of weak long range interactions, which is very difficult to design {\em de novo}. The alternative approach, that bypasses the requirement to consider the complexity of protein core, is to use the modular approach, where the desired structural scaffold is defined by largely independent interactions between the building modules, which could be used to build complex polyhedra. DNA nanotechnology (Seeman (2004) for example) has used this approach to design complex DNA assembles which are defined by the complementarity of antiparallel DNA strands. Typically several DNA strands are used to construct such structures although it has been recently demonstrated by Ko\v car et al. (2016) that DNA-based polyhedra could also be assembled from a single chain. The straightforward complementarity of base pairs is crucial for the use of DNA for the design of complex structures. On the other hand, it was shown by Gradi\v sar and Jerala (2011) that polypeptides form coiled-coil dimers, where the specificity between the interacting polypeptide chains can also be designed based on the combination of hydrophobic and electrostatic interactions to build a set of orthogonal building modules. Therefore the concatenation of coiled-coil modules connected by flexible linkers into a single chain might also mimic the design principles of DNA nanostructures for the self-assembling complex nanostructures.

Indeed in $2013$ Gradi\v sar et al. (2013) presented a novel polypeptide self-assembly strategy for nanostructure design. The main success of their research is a construction of a polypeptide self-assembling tetrahedron by concatenating $12$ coiled-coil-forming segments separated by flexible peptide hinges in a prescribed order. To be more precise, a single polypeptide chain consisting of $12$ segments was routed through $6$ edges of the tetrahedron in such a way that every edge was traversed exactly twice. In this way $6$ coiled-coil dimers were created and interlocked into a stable tetrahedral structure. 

While the required mathematical support for the particular case of the tetrahedron was already given by Gradi\v sar et al. (2013), Klav\v zar and Rus (2013) presented the first mathematical model for the investigations that could follow. It was observed that a polyhedron $P$ which is composed from a single polymer chain can be naturally represented with a graph $G(P)$ of the polyhedron and since in the self-assembly process every edge of $G(P)$ corresponds to a coiled-coil dimer, exactly two segments are associated with every edge of $G(P)$. The authors have then shown that a polyhedral graph $P$ can be realized by interlocking pairs of polypeptide chains if its corresponding graph $G(P)$ contains a closed walk which traverses every edge exactly twice (double traces which will be more formally defined later). Beside that, no edge should immediately be succeeded by the same edge in the opposite direction and a vertex sequence $u \rightarrow v \rightarrow w$ can appear at most once in any direction ($u \rightarrow v \rightarrow w$ or $w \rightarrow v \rightarrow u$) in the double trace. Fijav\v z et al. (2014) observed that the model presented by Klav\v zar and Rus (2013) is the appropriate mathematical description for graphs already constructed from coiled-coil-forming segments, yet it is deficient in modeling graphs with either very small ($\leq 2$) or large ($\geq 6$) degree vertices. Since the goal is to realize such graphs with coiled-coil-forming segments in the future, strong and $d$-stable traces (to be defined later) were presented as a natural extension of the model by Klav\v zar and Rus (2013) and graphs admitting them were characterized.

Even before, self-assembly strategy of designed DNA into different polyhedra were reported: into tetrahedron by Goodman et al. (2004) and by He (2008), into cube by Chen and Seeman (1991) and by Zhang et al. (2009), into octahedron by Shih et al. (2004), into dodecahedron by Zimmermann (2008), and into icosahedron by Bhatia et al. (2009) and by Douglas et al. (2009). In order to expand a novel polypeptide self-assembly strategy to self-assembly of designed DNA and in order to have a self-assembly strategy of designed DNA into arbitrary polyhedron, the question when strong and $d$-stable traces which traverse every edge once in each direction was raised. While graphs which admit antiparallel strong traces were characterized by Fijav\v z et al. (2014), one of the open problems left there was a charecterization of graphs admitting antiparallel $d$-stable traces. An {\em antiparallel $d$-stable trace} is a closed walk which traverses every edge exactly once in each direction and for every vertex $v$, there is no subset $N$ of its neighbors, with $1 \leq |N| \leq d$, such that every time the walk enters $v$ from $N$, it also exits to a vertex in $N$. 

Let $T$ be a subtree of a graph $G$. Recall that an edge complement $G - E(T)$ of a subtree of $G$ is called a {\em co-tree}. Our main result, a characterization of graphs admitting antiparallel $d$-stable traces, can then be read as follows.

\begin{theorem}
\label{thm:main}
Let $d \geq 1$ be an integer. A graph $G$ admits an antiparallel $d$-stable trace if and only if $\delta(G) > d$ and $G$ has a spanning tree $T$ such that each component of the co-tree $G-E(T)$ has an even number of edges or contains a vertex $v$ of degree at least $2d + 2$.
\end{theorem}

Actually related topics were, motivated with embedding and single face embeddings of graphs in surfaces, studied a long time ago. In $1951$, Ore (1951) posed a problem, asking for a characterization of graphs that admit closed walk which traverses every edge exactly once in each direction and such that no edge is succeeded by the same edge in the opposite direction. The problem was partially solved by Troy (1966) and Eggleton and Skilton (1984), and completely solved almost $40$ years later by Thomassen (1990) as follows:

\begin{theorem}
\label{thm:1anti}
{\rm [Theorem~$3.3$ of Thomassen (1990)]}
A graph $G$ admits a closed walk which traverses every edge exactly once in each direction and such that no edge is succeeded by the same edge in the opposite direction if and only if $\delta(G) > 1$ and $G$ has a spanning tree $T$ such that each component of the co-tree $G-E(T)$ has an even number of edges or contains a vertex $v$ of degree $\geq 4$.
\end{theorem}

Note that Theorem~\ref{thm:1anti} is Theorem~\ref{thm:main} for $d = 1$. 

Independently, several generalizations of Theorem~\ref{thm:1anti} motivated by biomolecular computing and assembling graphs from strands of DNA were already introduced in the past. Ellis-Monaghan (2004) observed how Theorem~\ref{thm:1anti} characterizes graphs which may be constructed from a single strand of DNA. A generalization of Theorem~\ref{thm:1anti} different from ours presented by Fan and Zhu (1998) was used to characterize graphs which may be constructed from $m$ strands of DNA.

The paper is organized as follows. In the next section we first list basic concepts needed throughout the paper. Then we present double traces, $d$-repetitions, $d$-stable traces, strong traces, and some already known results about them. In Section~\ref{sec:spanning_trees} we prove some results about spanning trees for which also co-trees fulfill some additional requirements. Those results are later used in Section~\ref{sec:proof} where we prove Theorem~\ref{thm:main}. We conclude with some observations which follow from the theorem. Among others, an alternative visualization of graphs from Theorem~\ref{thm:main} is observed --- embedding graphs in pseudosurfaces, which in the past has received much less attention than embedding graphs in surfaces.  

%%%%%%%%%%%%%%%%%%
\section{Double traces}
\label{sec:double_traces}

All graphs considered in this paper will be finite. We denote the degree of a vertex $v$ by $d_G(v)$ or $d(v)$ for short if graph $G$ is clear from the context. The minimum and the maximum degree of $G$ are denoted with $\delta(G)$ and $\Delta(G)$, respectively. If $v$ is a vertex then $N(v)$ denotes a set of vertices adjacent to $v$, and $E(v)$ is the set of edges incident with $v$. 

A \emph{walk} in $G$ is an alternating sequence 
\begin{equation}
W=v_0 e_1 v_1 \ldots v_{\ell-1} e_\ell v_\ell,
\label{eq:walk}
\end{equation}
so that for each $i=1,\ldots,\ell$, $e_i$ is an edge between vertices $v_{i-1}$ and $v_i$. We say that $W$ \emph{passes through} or \emph{traverses} edges and vertices contained in the sequence~\eqref{eq:walk}. The length of a walk is the number of edges in the sequence, and we call $v_0$ and $v_\ell$ the \emph{endvertices} of $W$. A walk is \emph{closed} if its endvertices coincide. 

An \emph{Euler tour} in $G$ is a closed walk which traverses every edge of $G$ exactly once. $G$ is an \emph{Eulerian graph} if it admits an Euler tour. The fundamental Euler's theorem asserts that a (connected) graph $G$ is Eulerian if and only if all of its vertices are of even degree. A \emph{double trace} in $G$ is a closed walk which traverses every edge of $G$ exactly twice. By doubling every edge of graph, thus creating an Eulerian graph, and considering an Euler tour in this graph, Euler and later K\" onig (1936) observed that every connected graph $G$ has a double trace.

Let $W$ be a double trace of length $\ell$ in graph $G$, $v$ an arbitrary vertex in $G$, and $N \subseteq N(v)$. $W$ has an \emph{$N$-repetition at $v$} if the following implication holds:
\begin{equation}
\text{\emph{for every $i \in \{0,\ldots,\ell-1\}$: if $v=v_i$, then $v_{i+1} \in N$ if and only if $v_{i-1} \in N$.}}
\label{eq:repetition}
\end{equation}
In other words, a double trace $W$ has an $N$-repetition at $v$ if whenever $W$ visits $v$ coming from a vertex in $N$ it also returns to a vertex in $N$. Note that $v_1$ is the vertex immediately following $v_\ell$, since we treat a double trace as a closed walk taking indices in~\eqref{eq:repetition} modulo $\ell$.

An $N$-repetition (at $v$) is a \emph{$d$-repetition} if $|N|=d$, and a $d$-repetition is sometimes also called a repetition \emph{of order $d$}. An $N$-repetition at $v$ is \emph{trivial} if $N=\emptyset$ or $N=N(v)$. Clearly if $W$ has an $N$-repetition at $v$, then it also has an $(N(v)\setminus N)$-repetition at $v$, a property also called a \emph{symmetry of repetitions}.

Fijav\v z et al. (2014) called double traces without nontrivial repetitions of order $\leq d$ as \emph{$d$-stable traces} while a \emph{strong trace} was defined as a double trace without nontrivial repetitions. Here the same notations will be used.

The next fact easily follows from the definitions of the $d$-stable trace, the strong trace and the symmetry of repetitions. For $d = 2$ it was already observed by Fijav\v z et al. (2014).

\begin{fact}
\label{fac:strong-stable}
In a graph $G$ with $\delta(G) > d$ every strong trace is a $d$-stable trace. If also $\Delta(G) < 2d + 2$, then every $d$-stable trace is also a strong trace.
\end{fact}

Graphs admitting strong and $d$-stable traces were characterized by Fijav\v z et al. (2014) using a correspondence between strong traces and single face embeddings. 

\begin{theorem}
\label{thm:strong}
{\rm [Theorem~$2.4$ of Fijav\v z et al. (2014)]}
Every connected graph $G$ admits a strong trace.
\end{theorem}

\begin{proposition}
\label{prop:dstable}
{\rm [Proposition~$3.4$ of Fijav\v z et al. (2014)]}
Let $G$ be a connected graph. Then $G$ admits a $d$-stable trace if and only if $\delta(G)>d$.
\end{proposition}

In an arbitrary double trace $W$ of a graph $G$ every edge is traversed twice. If $W$ traverses an edge $e=uv$ in the same direction twice (either both times from $u$ to $v$ or both times from $v$ to $u$) then we call $e$ a {\em parallel edge}, otherwise $e$ is  an {\em antiparallel edge}. A double trace $W$ is a {\em parallel double trace} if every edge of $G$ is parallel and an {\em antiparallel double trace} if every edge of $G$ is antiparallel. The next theorem proven by Fijav\v z et al. (2014) characterizes graphs which admit antiparallel strong traces. It turns out that those are exactly the graphs which have a single face embedding in an orientable surface. 

\begin{theorem}
\label{thm:santi}
{\rm [Theorem~$4.1$ of Fijav\v z et al. (2014)]}
A graph $G$ admits an antiparallel strong trace if and only if $G$ has a spanning tree $T$ such that each component of the co-tree $G - E(T)$ has an even number of edges.
\end{theorem}

To make our proofs more transparent we conclude this section with the next concept similar to the {\em detachment of graph} presented by Nash-Williams (1979, 1985, 1987) in a series of papers.

\begin{definition}
Let $G$ be a graph, $v$ a vertex of degree $\geq 2$, and $(N_1, \ldots, N_k)$ a partition of $N(v)$. We obtain graph $G'$ from graph $G$ with {\em splitting procedure} in $v$ using $N_1, \ldots, N_k$ as follows. Replace vertex $v$ with $k$ new nonadjacent vertices $v_1, \ldots, v_k$ in $G'$. Add edges between $v_i$ and the vertices from $N_i$ for $1 \leq i \leq k$. See also Fig.~\ref{fig:construction3}.
\end{definition}

For any other details about double traces, graph embeddings, or other terms and concepts from graph theory not defined here we refer to the books by Fleischner (1990, 1991), by Mohar and Thomassen (2001), and by West (1996), respectively.

%%%%%%%%%%%%%%%%%%%%%%%%%%%%%%%%%%%%%
\section{Spanning trees, co-trees and deficiency}
\label{sec:spanning_trees}

A tree $T$ is a spanning tree if and only if $G - E(T)$ is a minimum co-tree. The number of edges in any spanning tree of $G$ is $|V(G)| - 1$, while the number of edges in any minimum co-tree of $G$ is equal to $|E(G)| - |V(G)| + 1$ and is called a {\em Betti number}, $\beta(G)$. Note that a co-tree is not necessarily connected. A component $C$ of a co-tree $G-E(T)$ is called an {\em even component} or an {\em odd component} if $C$ has an even or an odd number of edges, respectively. The {\em deficiency} of a spanning tree $T$ in $G$, denoted with $\xi(G, T)$, is defined as the number of odd components of the co-tree $G - E(T)$. The {\em deficiency of a graph $G$}, denoted with $\xi(G)$, is defined as the minimum tree deficiency over all spanning trees $T$ of $G$. Spanning trees that realize $\xi(G)$ are called {\em Xuong trees}. Spanning trees, co-trees, Betti number, and deficiency were used for new graph embeddings technique by Xu (1979).

Here we are interested in spanning trees of $G$ such that each odd component of the co-tree contains a high degree vertex $v$ ($d_G(v) \geq d$ for some integer $d$). Therefore we define {\em $d$-deficiency of a graph $G$}, denoted with $\xi(G,d)$, as the minimum deficiency over all spanning trees $T$ of $G$ such that each odd component of the co-tree $G-E(T)$ contains a vertex $v$ with $d_G(d) \geq d$. Note that for every graph $G$ and any integer $d$, $\xi(G) \leq \xi(G,d)$. 

For such a spanning tree $T$ in $G$ that each component of the co-tree $G - E(T)$ is even or contains a vertex of degree at least $d$, we also define {\em $d$-deficiency} of spanning tree $T$ in $G$, denoted with $\xi(G,d,T)$, as the number of odd component in the co-tree $G - E(T)$.

Since any vertex $v \in V(G)$ is contained in at most one (odd) component of the co-tree $G - E(T)$, the next definition makes sense: 

\begin{definition}
\label{def:odd_comp}
Let $G$ be a graph, $T$ a spanning tree of $G$, $C$ an odd component in the co-tree $G -E(T)$, and $v$ an arbitrary vertex in $C$. Let $v_1, \ldots, v_c$ be the neighbors of $v$ in $C$. $\mathcal{O}(v,T)$ is the set of odd components that we obtain from $C$ if we split $v$ into $c$ new vertices and pairwise adjacent them to $v_1, \ldots, v_c$. Analogously we denote the set of even components with $\mathcal{E}(v,T)$. 
\end{definition}

Observe that because $C$ is an odd component, the number of odd components in $\mathcal{O}(v,T)$ is also odd and therefore $|\mathcal{O}(v,T)| \geq 1$. In the rest of this section we prove some observations about spanning trees of graphs.

\begin{lemma}
\label{lemma:treeContractE}
Let $G$ be a connected graph, $U \subseteq V(G)$, $v$ an arbitrary vertex of degree $\geq 2$ in $G$, and $N_1, \ldots, N_k \subseteq N(v)$, $k \geq 2$, a partition of neighbors of $v$. Let graph $G'$ be obtained from $G$ with splitting procedure in $v$ using $N_1, \ldots, N_k$. Let $v_1, \ldots, v_k$ be new vertices in $G'$. If $G'$ has a spanning tree $T'$ such that each component of the co-tree $G' - E(T')$ is even or contains a vertex from $U$ or at least one of new vertices $v_1, \ldots, v_k$, then there exists a spanning tree $T$ in $G$ such that each component of the co-tree $G - E(T)$ is even or contains a vertex from $U$ or contains $v$.  
\end{lemma}

\proof
Assume that $G$ and $G'$ are two connected graphs such that $G$ can be obtained from $G'$ by identifying $k \geq 2$ vertices $v_1, \ldots, v_k$ that have disjoint neighborhoods with an arbitrary vertex $v$ of degree at least $2$. In addition assume that $G'$ has a spanning tree $T'$ such that  each component of the co-tree $G' - E(T')$ is even or contains a vertex from $U \subseteq V(G) \setminus \{v_1, \ldots, v_k\}$ or contains at least one of $v_1, \ldots, v_k$. 

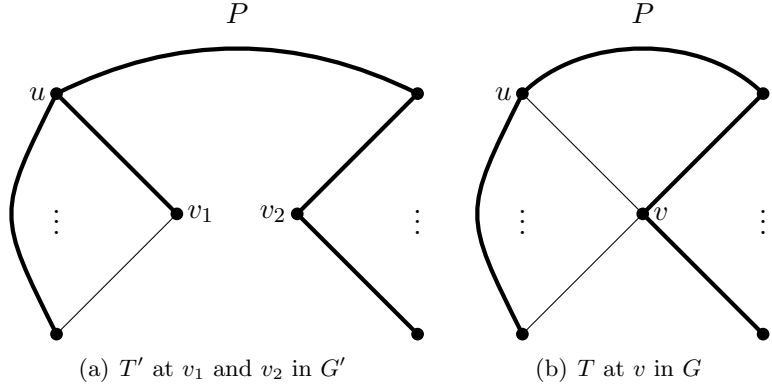
\begin{figure}[ht!]
\begin{center}
\subfigure[$T'$ at $v_1$ and $v_2$ in $G'$]
{
\begin{tikzpicture}[scale=0.8,style=thick]
\fill (2,2) circle (3pt) node[right]{$v_1$};
\fill (4,2) circle (3pt) node[left]{$v_2$};
\fill (0,0) circle (3pt);
\fill (0,2) node{$\vdots$};
\fill (0,4) circle (3pt) node[left]{$u$};
\fill (6,0) circle (3pt);
\fill (6,2) node{$\vdots$};
\fill (6,4) circle (3pt);
\draw[thin] (2,2)--(0,0);
\draw[ultra thick] (2,2)--(0,4);
\draw[ultra thick] (4,2)--(6,0);
\draw[ultra thick] (4,2)--(6,4);
\fill (3,5) node[above]{$P$};
\draw[ultra thick] (0,4)  .. controls (2,5) and (4,5)  .. (6,4);
\draw[ultra thick] (0,4)  .. controls (-1,2) and (-1,2)  .. (0,0);
\end{tikzpicture}
}
\subfigure[$T$ at $v$ in $G$]
{
\begin{tikzpicture}[scale=0.8,style=thick]
\fill (2,2) circle (3pt) node[right]{$v$};
\fill (0,0) circle (3pt);
\fill (0,2) node{$\vdots$};
\fill (0,4) circle (3pt) node[left]{$u$};
\fill (4,0) circle (3pt);
\fill (4,2) node{$\vdots$};
\fill (4,4) circle (3pt);
\draw[thin] (2,2)--(0,0);
\draw[thin] (2,2)--(0,4);
\draw[ultra thick] (2,2)--(4,0);
\draw[ultra thick] (2,2)--(4,4);
\fill (2,5) node[above]{$P$};
\draw[ultra thick] (0,4)  .. controls (1,5) and (3,5)  .. (4,4);
\draw[ultra thick] (0,4)  .. controls (-1,2) and (-1,2)  .. (0,0);
\end{tikzpicture}
}
\end{center}
\caption{Construction of a spanning tree $T$ in $G$ from a spanning tree $T'$ in $G'$ from the proof of Lemma~\ref{lemma:treeContractE}. Edges contained in trees $T'$ and $T$ are drawn thick.}
\label{fig:construction1}
\end{figure}

We proceed by induction on $k$. Let $k = 2$ and let $v_1$ and $v_2$ be two vertices that have disjoint neighborhoods $N(v_1)$ and $N(v_2)$ that we have to identify into $v$ in order to obtain $G$ from $G'$. Obtain subgraph $T''$ in $G$ from $T'$ as follows. Let $e' = xy$ be an arbitrary edge in $T'$. If $x,y \notin \{v_1,v_2\}$ we put $e'$ into $T''$. If $x = v_1$ then replace $e'$ with $vy$. Analogously we replace edge where $x = v_2$, or $y = v_1$, or $y = v_2$ with $vy$, $xv$, and $xv$, respectively (if they appear in $T'$). Since $T'$ is a spanning tree of $G'$, there exists a unique $v_1,v_2$-path $P$ in $T'$. Let $u$ be the neighbor of $v_1$ in $P$ and $e = uv_1$, see Fig.~\ref{fig:construction1} (a). 

We claim that $T = T'' - e$ is a spanning tree of $G$ and each component of the co-tree $G - E(T)$ is even or contains a vertex from $U$ or contains $v$. If $T$ is not a spanning subgraph of $G$ then there exists a vertex $x \in G \setminus T$. Since $v$ is incident with at least one edge from $T$ by construction, $x \neq v$. If $x \in V(G) \setminus \{v\}$ is not incident with any edge from $T$ it follows that $x$ is not incident with any edge from $T'$ in $G'$, a contradiction. Therefore $T$ is a spanning subgraph in $G$. Let now $x, y \neq v$ be two arbitrary vertices in $G$. Since $T'$ is a spanning tree of $G'$, there exists a unique $x,y$-path $Q'$ in $T'$. If $Q'$ does not contain an edge $e$ then it is clear from the construction of $T$ that $Q'$ is also contained in $T$. Assume now that $e \in Q'$. Without loss of generality we can write $Q' = x-Q_1'-v_1-e-u-Q_2'-y$. Denote the part of the path $P$ that goes from $v_2$ to $u$ with $P'$. Since $u$ is an endvertex of $Q_2'$ and $u$ is an endvertex of $P'$, $Q_2'$ and $P'$ are not disjoint. We can then find a $x,y$-path in $x-Q_1'-v-P'-u-Q_2'-y$. Analogously if $x = v$ or $y = v$. Therefore $T$ is connected. Let now $C$ be a cycle in $T$. If $C$ does not contain the vertex $v$ it follows from the construction of $T$ that $C$ is a cycle in $T'$, which is absurd. Therefore $C$ has to contain $v$. Let $x$ be an arbitrary vertex in $C$ different from $v$. Without loss of generality we can write $C = v-C_1-x-C_2-v$, where $C_1$ and $C_2$ are two vertex disjoint $x,v$-paths. If $C_1$ and $C_2$ are both disjoint $v,v_1$-paths in $T'$ or both disjoint $x,v_2$-paths in $T'$ we get a contradiction with the fact that $T'$ is a spanning tree. Therefore without loss of generality we can assume that $C_1$ is a $x,v_1$-path in $T'$ while $C_2$ is a $x,v_2$-paths in $T'$. Since $T'$ is a spanning tree it follows that for any two vertices there exists a unique path between them in $T'$. Therefore $P = v_1-C_1-x-C_2-v_2$ and $C_1$ starts with $e$. Since $T = T'' - e$, $C_1$ can not be contained in $T$. It follows that $T$ is a spanning tree of $G$. 

Assume that $O$ is an odd component of the co-tree $G - E(T)$ which does not contain a vertex from $U$ or $v$. It follows that $O$ is an odd component of the co-tree $G' - E(T')$ which does not contain a vertex from $U$ or any of the vertices $v_1$ and $v_2$ which is absurd. We have thus proved that if $G$ is obtained from $G'$ by identifying exactly two vertices into vertex $v$, then $G$ has a spanning tree $T$ such that each component of the co-tree $G - E(T)$ is even or contains a vertex from $U$ or contains $v$. 

Let next $k > 2$, let induction hypothesis be true for any $l < k$, and let $v_1, \ldots, v_k$ be vertices that have disjoint neighborhoods $N(v_1), \ldots, N(v_k)$, that we have to identify into $v$ in order to obtain $G$ from $G'$. Construct $G''$ from $G'$ by identifying vertices $v_1, \ldots, v_{k-1}$ into a new vertex $v''$. By induction $G''$ has a spanning tree $T''$ such that each component of the co-tree $G'' - E(T'')$ is even or contains a vertex from $U$ or contains $v''$. Identify vertices $v_k$ and $v''$ into a vertex $v$ in $G''$ to obtain the graph $G$. Now by induction also $G$ has a spanning tree $T$ such that each component of the co-tree $G - E(T)$ is even or contains a vertex from $U$ or contains $v$. 
\qed

Note that in Lemma~\ref{lemma:treeContractE}, $U$ can be the empty set.

\begin{lemma}
\label{lemma:treeSplittV}
Let $G$ be a connected graph, $T$ a spanning tree of $G$, and $v$ an arbitrary vertex in $G$ such that $v$ is contained in an odd component of the co-tree $G - E(T)$. Let $\mathcal{G}$ be the family of all graphs obtained from $G$ with splitting procedure in $v$ using two subsets of $N(v)$ of cardinality $\left\lceil\frac{d(v)}{2}\right\rceil$ and $\left\lfloor\frac{d(v)}{2}\right\rfloor$. There exists a graph $G' \in \mathcal{G}$ such that $G'$ has a spanning tree $T'$ and $\xi(G', T') < \xi(G, T)$. 
\end{lemma}

\proof
Let $G$ be a connected graph, $T$ a spanning tree of $G$, and $v$ an arbitrary vertex in $G$ such that $v$ is contained in an odd component of the co-tree $G - E(T)$. Since $T$ is a spanning tree of $G$ there exists an edge $e_y = yv \in E(T)$ and because $v$ is contained in an odd component of the co-tree $G-E(T)$, there also exists an edge $e_w = wv \in E(G) \setminus E(T)$. Clearly $e_y$ and $e_w$ are two different edges, therefore $d_G(v) \geq 2$. 

\begin{figure}[ht!]
\begin{center}
\subfigure[$T$ at $v$ in $G$]
{
\begin{tikzpicture}[scale=0.8,style=thick]
\fill (2,2) circle (3pt) node[right]{$v$};
\fill (0,0) circle (3pt);
\fill (0,2) node{$\vdots$};
\fill (0,4) circle (3pt) node[left]{$u$};
\fill (4,0) circle (3pt);
\fill (4,2) node{$\vdots$};
\fill (4,4) circle (3pt) node[right]{$w$};
\draw[thin] (2,2)--(0,0);
\draw[ultra thick] (2,2)--(0,4);
\draw[ultra thick] (2,2)--(4,0);
\draw[thin] (2,2)--(4,4);
\fill (1,3) node[above right]{$e_u$};
\fill (3,3) node[above left]{$e_w$};
\fill (2,5) node[above]{$P$};
\draw[ultra thick] (0,4)  .. controls (1,5) and (3,5)  .. (4,4);
\draw[ultra thick] (0,4)  .. controls (-1,2) and (-1,2)  .. (0,0);
\end{tikzpicture}
}
\subfigure[$T'$ at $v'$ and $v''$ in $G'$]
{
\begin{tikzpicture}[scale=0.8,style=thick]
\fill (2,2) circle (3pt) node[right]{$v'$};
\fill (4,2) circle (3pt) node[left]{$v''$};
\fill (0,0) circle (3pt);
\fill (0,2) node{$\vdots$};
\fill (0,4) circle (3pt) node[left]{$u$};
\fill (6,0) circle (3pt);
\fill (6,2) node{$\vdots$};
\fill (6,4) circle (3pt) node[right]{$w$};
\draw[thin] (2,2)--(0,0);
\draw[ultra thick] (2,2)--(0,4);
\draw[ultra thick] (4,2)--(6,0);
\draw[ultra thick] (4,2)--(6,4);
\fill (1,3) node[above right]{$e_u$};
\fill (5,3) node[above left]{$e_w$};
\fill (3,5) node[above]{$P$};
\draw[ultra thick] (0,4)  .. controls (2,5) and (4,5)  .. (6,4);
\draw[ultra thick] (0,4)  .. controls (-1,2) and (-1,2)  .. (0,0);
\end{tikzpicture}
}
\end{center}
\caption{Construction of a spanning tree $T'$ in $G'$ from a spanning tree $T$ in $G$ from the proof of Lemma~\ref{lemma:treeSplittV}. Edges contained in trees $T$ and $T'$ are drawn thick.}
\label{fig:construction2}
\end{figure}
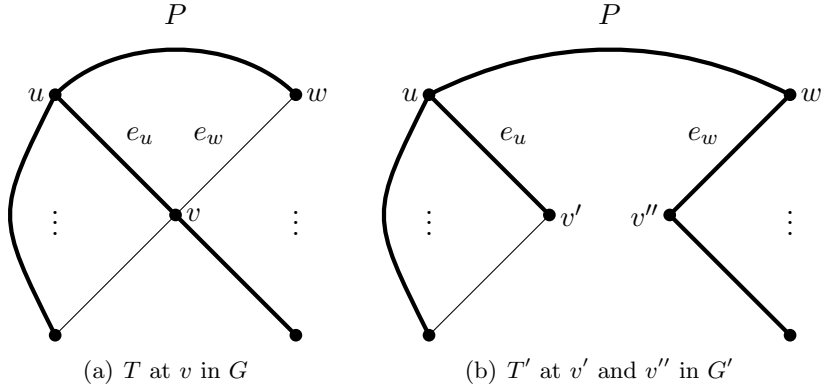

Because $T$ is a spanning tree of $G$, there exists a unique $v,w$-path in $T$. Denote the vertex adjacent to $v$ in this path with $u$ and the part of this path between $u$ and $w$ with $P$. Let $N$ be a subset of $N(v)$ with $\left\lceil\frac{d(v)}{2}\right\rceil$ vertices such that $u \in N$ while $w \notin N$. Obtain a graph $G'$ from $G$ with the splitting procedure in $v$ using $N$ and $N(v) \setminus N$. Denote new vertices which are connected to vertices from $N$ and $N(v) \setminus N$ with $v'$ and $v''$, respectively. Construct the subgraph $T'$ in $G'$ from $T$ as follows. Let $e = xy$ be an arbitrary edge in $T$. If $x, y \neq v$ put $e$ into $T'$. If $x = v$ and $y \in N$ then replace $e$ with $v'y$. Analogously we replace edges where $x = v$ and $y \in N(v) \setminus N$, or $y = v$ and $x \in N$, or $y = v$ and $x \in N(v) \setminus N$ with $v''y$, $xv'$, or $xv''$, respectively. Finally put $e_w$ into $T'$. See Fig.~\ref{fig:construction2}. Parallel as in the proof of Lemma~\ref{lemma:treeContractE} we can prove that $T'$ is a spanning tree of $G'$. 

Because $u$ and $w$ were arbitrary neighbors of $v$ in $G$, such that $uv \in E(T)$ and $wv \in E(G) \setminus E(T)$, every graph obtained from $G$ with splitting procedure in $v$ using $u$ and $w$ with above described properties has a spanning tree $T'$ (obtained from $T$ by adding edge $wv$) and is therefore connected. We claim that for at least one of those graphs $G'$, $T'$ is such a spanning tree that $\xi(G', T') < \xi(G, T)$.

Let from now on $C$ be the odd component of the co-tree $G - E(T)$ which contains $v$. Note first that because $v$ lies in exactly one component of the co-tree $G - E(T)$, every component $C'$ in the co-tree $G-E(T)$ different from $C$ will also be in the co-tree $G'-E(T')$. Note next that components which will be formed from $C$ in $G' - E(T')$ will be constructed from components in $\mathcal{E}(v,T)$ and $\mathcal{O}(v,T)$ from one of which one of the edges incident with $v$ will be removed. Finally, it is clear from Definition~\ref{def:odd_comp} that every component $O \in \mathcal{O}(v,T)$ and $E \in \mathcal{E}(v,T)$ remains connected if we remove an edge incident with $v$ from $O$. Note also that $O$ and $E$ can have more than one edge incident with $v$.

We consider two cases:

\vspace{5pt}
\noindent
\textbf{Case 1:} $|\mathcal{O}(v,T)| \leq \left\lceil\frac{d(v)}{2}\right\rceil$

\vspace{5pt}
\noindent
Let $e = vw \in O_w$ for some $O_w \in \mathcal{O}(v,T)$ and $f = vu \in T$ be two edges incident with $v$ in $G$. Use them as $e_w$ and $e_y$ in the at the beginning of the proof described construction of $G'$, respectively. Since $|\mathcal{O}(v,T)| \leq \left\lceil\frac{d(v)}{2}\right\rceil$ we can connect $v''$ to $w$ and for every from $O_w$ different component $O \in \mathcal{O}(v,T)$ to at least one vertex from $O$, when constructing $G'$. When $e_w$ is put in $T'$, the odd component $O_w$ becomes even. In $\mathcal{O}(v,T)$ therefore remains an even number (possibly $0$) of odd components which are all incident with $v''$ in $G'$, therefore they still belong to one component of $G' - E(T')$, which is now even. Components from $\mathcal{E}(v,T)$ also remain even. Therefore $T'$ is a spanning tree of $G'$ such that $\xi(G', T') < \xi(G, T)$.

\vspace{5pt}
\noindent
\textbf{Case 2:} $|\mathcal{O}(v,T)| > \left\lceil\frac{d(v)}{2}\right\rceil \geq 1$

\vspace{5pt}
\noindent
Since $|\mathcal{O}(v,T)| > 1$, $|\mathcal{O}(v,T)|$ is odd, and there also exists an edge in $T$ incident with $v$, it follows that $d(v) \geq 4$. We will use $uv$ to denote an edge incident with $v$ from $T$. We first describe how $G'$ is obtained from $G$ and $T'$ from $T$ such that each component of the co-tree $G'-E(T')$ is even if at least one of the next conditions is fulfilled:

\begin{itemize}
\item[$(i)$]
There exists an odd component $O_1 \in \mathcal{O}(v,T)$ with at least two edges incident with $v$. 
\item[$(ii)$]
$d(v) \neq 0 \pmod{4}$.
\item[$(iii)$]
Spanning tree $T$ and $v$ have more than one incident edge.
\item[$(iv)$]
$|\mathcal{E}(v,T)| > 0$.
\end{itemize}

In $(i)$ denote edges from $O_1$ with $e_1 = vw_1$ and $e_2 = vw_2$. Denote an edge incident with $v \in T$ with $f = vu \in T$. Since $|\mathcal{O}(v,T)| > 1$, there exists at least one more vertex $w_3$ such that $e_3 = vw_3 \in O_2 \in \mathcal{O}(v,T)$, where $O_1 \neq O_2$. Use $e_3$ and $f$ as $e_w$ and $e_y$ in the at the beginning of the proof described construction of $G'$, respectively. Since $\left\lceil\frac{d(v)}{2}\right\rceil, \left\lfloor\frac{d(v)}{2}\right\rfloor \geq 2$, we can connect $v'$ to $u$ and $w_1$, and connect $v''$ to $w_3$ and $w_2$. When $e_3$ is put in $T'$, odd component $O_2$ becomes an even component and since $w_1$ is connected to $v'$ and $w_2$ is connected to $v''$, $v'$ and $v''$ are in the same component of the co-tree $G' - E(T')$ . In $\mathcal{O}(v,T)$ we now have an even number of odd components. Therefore $T'$ is a spanning tree of $G'$ such that $\xi(G', T') < \xi(G, T)$. 

For $(ii)-(iv)$ the parity of $d(v)$, additional edge incident with $v$ from $T$, or an edge incident with $v$ from a component in $\mathcal{E}(v,T)$ is used to achieve that $v''$ is adjacent to an odd number of different components from $\mathcal{O}(v,T)$, while $v'$ is connected to $u$ via edge replacing $vu$ from $T$ and the rest of the odd components (an even number of them). Note that because components from $\mathcal{O}(v,T)$ and $\mathcal{E}(v,T)$ can have more than one edge incident with $v$, $v'$ and $v''$ can both be adjacent to the same component from $\mathcal{O}(v,T)$ or $\mathcal{E}(v,T)$. Note also that to achieve above described partition, we also allow that $d(v') = \left\lfloor\frac{d(v)}{2}\right\rfloor$ while $d(v'') = \left\lceil\frac{d(v)}{2}\right\rceil$. Let $O$ be one of those components which will be adjacent to $v''$ and $e = vw$ an edge in it. Use $e$ as $e_w$ in the at the beginning of the proof described construction. Therefore $O$ will then become an even component and $v'$ and $v''$ will be both connected with even number of odd components from $\mathcal{O}(v,T)$.

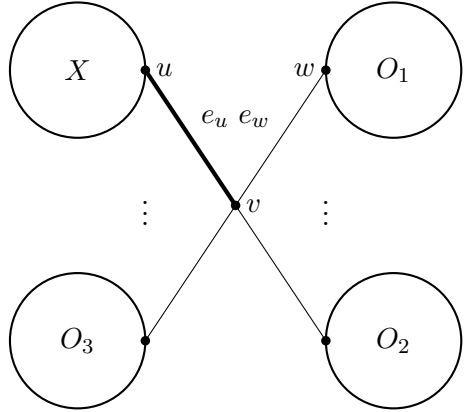
\begin{figure}[ht!]
\begin{center}
\begin{tikzpicture}[scale=0.6,style=thick]
\fill (2,3) circle (3pt) node[right]{$v$};
\fill (0,0) circle (3pt);
\fill (0,3) node{$\vdots$};
\fill (0,6) circle (3pt) node[right]{$u$};
\fill (4,0) circle (3pt);
\fill (4,3) node{$\vdots$};
\fill (4,6) circle (3pt) node[left]{$w$};
\draw (5.5,6) circle [radius=1.5] node{$O_1$};
\draw (5.5,0) circle [radius=1.5] node{$O_2$};
\draw (-1.5,6) circle [radius=1.5] node{$X$};
\draw (-1.5,0) circle [radius=1.5] node{$O_3$};
\draw[thin] (2,3)--(0,0);
\draw[ultra thick] (2,3)--(0,6);
\draw[thin] (2,3)--(4,0);
\draw[thin] (2,3)--(4,6);
\fill (1,4.5) node[above right]{$e_u$};
\fill (3,4.5) node[above left]{$e_w$};
\end{tikzpicture}
\end{center}
\caption{Graph $G$ and a spanning tree $T$ at vertex $v$ if none of the conditions from Case $2$ in the proof of Lemma~\ref{lemma:treeSplittV} is fulfilled. Note that because $T$ is a spanning tree, there exist $u,v$-paths in $T$ not shown in the figure, but there are no paths completely contained in $G-E(T)$ between different odd components. Edge contained in spanning tree $T$ is drawn thick.}
\label{fig:non_condition}
\end{figure}

Assume now that none of the above conditions is fulfilled. Then $|\mathcal{O}(v,T)| = d_G(v) - 1$ and we can without loss of generality assume that $G$ looks at $v$ as it is shown at Fig.~\ref{fig:non_condition}, where $O_1$, $O_2$, and $O_3$ at Fig.~\ref{fig:non_condition} are (together with edges incident with $v$) odd components from $\mathcal{O}(v,T)$. Note that there can be more than just $3$ odd components. Note also that $X$ which contains vertex $u$, can be an even (possibly empty) component of the co-tree $G-E(T)$ or an odd component disjoint with $C$ in the co-tree $G-E(T)$. 

\begin{figure}[ht!]
\begin{center}
\begin{tabular}{c c}
\subfigure[$T$ at $v$ in $G$]
{
\begin{tikzpicture}[scale=0.6,style=thick]
\fill (2,3) circle (3pt) node[right]{$v$};
\fill (0,0) circle (3pt)  node[right]{$v_3$};
\fill (0,3) node{$\vdots$};
\fill (0,6) circle (3pt) node[right]{$u$};
\fill (4,0) circle (3pt) node[left]{$v_2$};
\fill (4,3) node{$\vdots$};
\fill (4,6) circle (3pt) node[left]{$w$};
\fill (2,6.75) circle (3pt) node[below]{$p$};
\draw (5.5,6) circle [radius=1.5] node{$O_1$};
\draw (5.5,0) circle [radius=1.5] node{$O_2$};
\draw (-1.5,6) circle [radius=1.5] node{$X$};
\draw (-1.5,0) circle [radius=1.5] node{$O_3$};
\draw[thin] (2,3)--(0,0);
\draw[ultra thick] (2,3)--(0,6);
\draw[thin] (2,3)--(4,0);
\draw[thin] (2,3)--(4,6);
\fill (2,7) node[above]{$P_1$};
\draw[ultra thick] (0,6)  .. controls (1,7) and (3,7)  .. (4,6);
\fill (0.6,2.7) node[right]{$P_2$};
\draw[ultra thick, dashed] (0,6)  .. controls (1,3) and (3,1)  .. (4,0);
\fill (1,4.5) node[above right]{$e_u$};
\fill (3,4.5) node[above left]{$e_w$};
\end{tikzpicture}
}

&

\subfigure[$T'$ at $v'$ and $v''$ in $G'$]
{
\begin{tikzpicture}[scale=0.6,style=thick]
\fill (2,3) circle (3pt) node[right]{$v'$};
\fill (4,3) circle (3pt) node[left]{$v''$};
\fill (0,0) circle (3pt)  node[right]{$v_3$};
\fill (0,3) node{$\vdots$};
\fill (0,6) circle (3pt) node[right]{$u$};
\fill (2,7) circle (3pt) node[below]{$p$};
\fill (6,0) circle (3pt) node[left]{$v_2$};
\fill (6,3) node{$\vdots$};
\fill (6,6) circle (3pt) node[left]{$w$};
\draw (7.5,6) circle [radius=1.5] node{$O_1$};
\draw (7.5,0) circle [radius=1.5] node{$O_2$};
\draw (-1.5,6) circle [radius=1.5] node{$X$};
\draw (-1.5,0) circle [radius=1.5] node{$O_3$};
\draw[thin] (0,6)--(2,7);
\draw[ultra thick] (2,3)--(0,0);
\draw[thin] (2,3)--(0,6);
\draw[ultra thick] (4,3)--(6,0);
\draw[ultra thick] (4,3)--(6,6);
\fill (3,7) node[above]{$P_1$};
\draw[ultra thick] (2,7)  .. controls (4,7)  .. (6,6);
\fill (3,1.8) node[below]{$P_2$};
\draw[ultra thick, dashed] (0,6)  .. controls (0.5,2.5) and (5.5,0.5)  .. (6,0);
\fill (1,4.5) node[above right]{$e_u$};
\fill (5,4.5) node[above left]{$e_w$};
\end{tikzpicture}
}
\end{tabular}
\end{center}
\caption{First subcase in the proof of Lemma~\ref{lemma:treeSplittV}. There also exists a path in $T$ containing $v_3$ and only one of vertices $u$, $w$, or $v_2$ not shown in the Figure. Edges contained in trees $T$ and $T'$ are drawn thick.}
\label{fig:subcase1}
\end{figure}
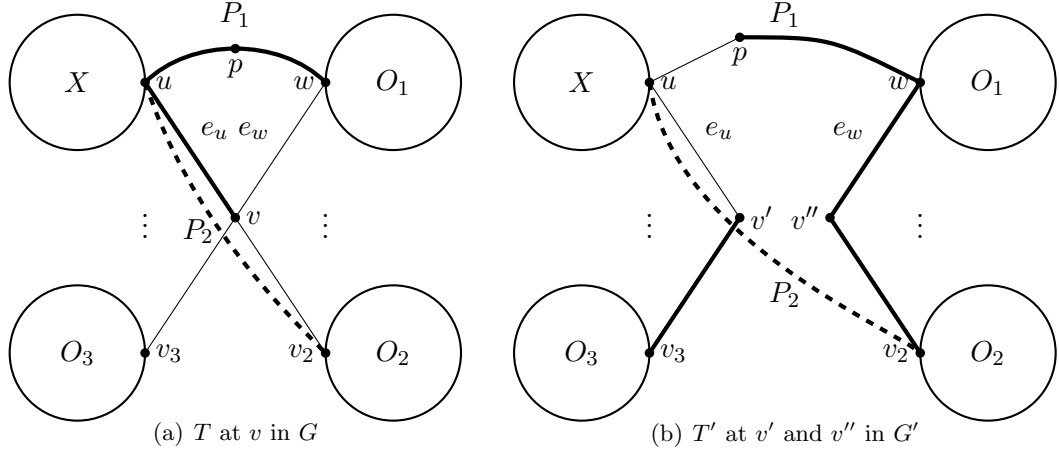

We consider two subcases. In the first subcase there exist such two odd components $O_1$ and $O_2$ in $\mathcal{O}(v,T)$, that there exists two disjoint paths $P_1$ and $P_2$ between them and vertex $u$. Connect $O_1$ and $O_2$ with $v''$ at vertices $w$ and $v_2$ in $G'$, respectively. Denote another component from $\mathcal{O}(v,T)$ which is in $G'$ connected to $v'$ with $O_3$ and a vertex from $O_3$ adjacent to $v$ with $v_3$. Denote the vertex adjacent to $u$ in $P_1$ with $p$ (note that $p$ can be equal to $w$). Construct spanning tree $T'$ from $T$ as described in previous subcases. Then add edges $v''v_2$ and $v'v_3$ to $T'$ and remove edges $v'u$ and $up$. See Fig.~\ref{fig:subcase1}.

It is not difficult to see that $T'$ is a spanning tree of $G'$. Component $O_2$ and $O_3$ are now even and are disjoint with every other odd component of $G' - E(T')$. Vertex $v''$ is still adjacent to an even number of odd components from $\mathcal{O}(v,T)$. Vertex $v'$ is still adjacent to an even number of odd components from $\mathcal{O}(v,T)$ that are different from $O_3$. Also $O_1$ is now even component of $G' - E(T')$, or is connected to $X$. The rest of the odd component $C$ has now common vertex $u$ with $X$. If $X$ was some from $C$ different odd component of $G-E(T)$, we connected $C$ with some other odd component of $G-E(T)$ in $T'$ and no matter if component containing $X$ is odd or even it follows that $\xi(G', T') < \xi(G, T)$. We have therefore connected two odd components into one odd component or join one even and one odd component into one even component (the parities are valid after we already consider removed edges).

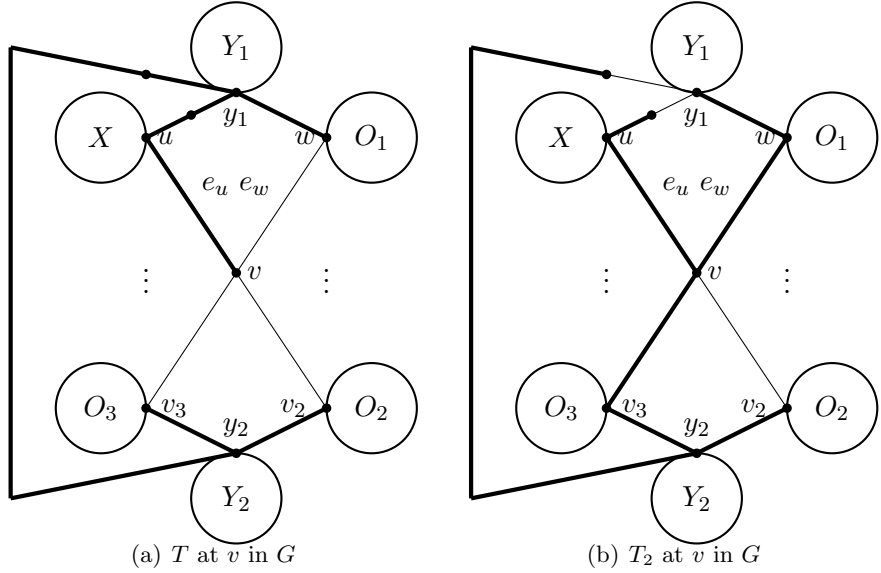
\begin{figure}[ht!]
\begin{center}
\begin{tabular}{c c}
\subfigure[$T$ at $v$ in $G$]
{
\begin{tikzpicture}[scale=0.6,style=thick]
\fill (2,-1) circle (3pt);
\fill (2,7) circle (3pt);
\fill (2,0) node[below]{$y_2$};
\fill (2,6) node[above]{$y_1$};
\fill (2,3) circle (3pt) node[right]{$v$};
\fill (0,0) circle (3pt);
\fill (0.1,0) node[right]{$v_3$};
\fill (0,3) node{$\vdots$};
\fill (0,6) circle (3pt) node[right]{$u$};
\fill (4,0) circle (3pt);
\fill (3.8,0) node[left]{$v_2$};
\fill (4,3) node{$\vdots$};
\fill (4,6) circle (3pt) node[left]{$w$};
\draw (5,6) circle [radius=1] node{$O_1$};
\draw (5,0) circle [radius=1] node{$O_2$};
\draw (-1,6) circle [radius=1] node{$X$};
\draw (2,-2) circle [radius=1] node{$Y_2$};
\draw (2,8) circle [radius=1] node{$Y_1$};
\draw (-1,0) circle [radius=1] node{$O_3$};
\draw[thin] (2,3)--(0,0);
\draw[ultra thick] (2,3)--(0,6);
\draw[thin] (2,3)--(4,0);
\draw[thin] (2,3)--(4,6);
\draw[ultra thick] (2,-1)--(0,0);
\draw[ultra thick] (2,-1)--(4,0);
\draw[ultra thick] (2,7)--(0,6);
\draw[ultra thick] (2,7)--(4,6);
\draw[ultra thick] (2,-1)--(-3,-2);
\draw[ultra thick] (-3,-2)--(-3,8);
\draw[ultra thick] (2,7)--(-3,8);
\fill (0,7.4) circle (3pt);
\fill (1,6.5) circle (3pt);
\fill (1,4.5) node[above right]{$e_u$};
\fill (3,4.5) node[above left]{$e_w$};
\end{tikzpicture}
}

&

\subfigure[$T_2$ at $v$ in $G$]
{
\begin{tikzpicture}[scale=0.6,style=thick]
\fill (2,-1) circle (3pt);
\fill (2,7) circle (3pt);
\fill (2,0) node[below]{$y_2$};
\fill (2,6) node[above]{$y_1$};
\fill (2,3) circle (3pt) node[right]{$v$};
\fill (0,0) circle (3pt);
\fill (0.1,0) node[right]{$v_3$};
\fill (0,3) node{$\vdots$};
\fill (0,6) circle (3pt) node[right]{$u$};
\fill (4,0) circle (3pt);
\fill (3.8,0) node[left]{$v_2$};
\fill (4,3) node{$\vdots$};
\fill (4,6) circle (3pt) node[left]{$w$};
\draw (5,6) circle [radius=1] node{$O_1$};
\draw (5,0) circle [radius=1] node{$O_2$};
\draw (-1,6) circle [radius=1] node{$X$};
\draw (2,-2) circle [radius=1] node{$Y_2$};
\draw (2,8) circle [radius=1] node{$Y_1$};
\draw (-1,0) circle [radius=1] node{$O_3$};
\draw[ultra thick] (2,3)--(0,0);
\draw[ultra thick] (2,3)--(0,6);
\draw[thin] (2,3)--(4,0);
\draw[ultra thick] (2,3)--(4,6);
\draw[ultra thick] (2,-1)--(0,0);
\draw[ultra thick] (2,-1)--(4,0);
\draw[thin] (2,7)--(1,6.5);
\draw[ultra thick] (1,6.5)--(0,6);
\draw[ultra thick] (2,7)--(4,6);
\draw[ultra thick] (2,-1)--(-3,-2);
\draw[ultra thick] (-3,-2)--(-3,8);
\draw[thin] (2,7)--(0,7.4);
\draw[ultra thick] (0,7.4)--(-3,8);
\fill (0,7.4) circle (3pt);
\fill (1,6.5) circle (3pt);
\fill (1,4.5) node[above right]{$e_u$};
\fill (3,4.5) node[above left]{$e_w$};
\end{tikzpicture}
}
\end{tabular}
\end{center}
\caption{Second subcase in the proof of Lemma~\ref{lemma:treeSplittV}. Edges contained in trees $T$ and $T_2$ are drawn thick.}
\label{fig:subcase2}
\end{figure}

In the second subcase two odd components in $\mathcal{O}(v,T)$ which would be connected to $u$ with two disjoint paths do not exists. Without loss of generality we can assume that $T$ looks like is shown at Fig.~\ref{fig:subcase2} (a). It is not difficult to see that a spanning tree $T_2$ at Fig.~\ref{fig:subcase2} (b) fulfills condition $(iii)$ and that $\xi(G, T_2) \leq \xi(G, T)$. Analogously when $y_1 = y_2$ or $y_1 = w$ (and $Y_1 = Y_2$ or $Y_1 = O_1$).

We have thus constructed graph $G'$ which has a spanning tree $T'$ such that $\xi(G', T') < \xi(G, T)$ and since $\xi(G,T)$ was not in ahead prescribed, above described construction is always possible. 
\qed

Because vertex $v$ from Lemma~\ref{lemma:treeSplittV} is arbitrary and since during the constructing from the the proof of Lemma~\ref{lemma:treeSplittV}, $G \setminus (v \cup N(v))$ remains equivalent, next lemma easily follows if we assume that vertex $v$ is of degree $\geq d$:

\begin{lemma}
\label{lemma:treeSplittVd}
Let $G$ be a connected graph, let $T$ be such a spanning tree of $G$ that each component of the co-tree $G - E(T)$ is even or contains a vertex of degree at least $d$, and $v$ an arbitrary vertex of degree at least $d$ in $G$ which is contained in an odd component of the co-tree $G - E(T)$. Let $\mathcal{G}$ be the family of all graphs obtained from $G$ with splitting procedure in $v$ using two subsets of $N(v)$ of cardinality $\left\lceil\frac{d(v)}{2}\right\rceil$ and $\left\lfloor\frac{d(v)}{2}\right\rfloor$. There exists a graph $G' \in \mathcal{G}$ such that $G'$ is connected and it has a spanning tree $T'$ such that each component of the co-tree $G' - E(T')$ is even or contains a vertex of degree at least $d$ and $\xi(G',d,T') < \xi(G,d,T)$. 
\end{lemma}

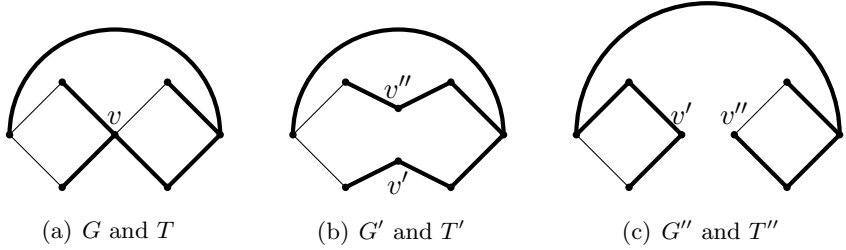
\begin{figure}[ht!]
\begin{center}
\begin{tabular}{c c c}
\subfigure[$G$ and $T$]
{
\begin{tikzpicture}[scale=0.7,style=thick]
\draw[thin] (0,1)--(1,2);
\draw[thin] (0,1)--(1,0);
\draw[ultra thick] (1,0)--(2,1);
\draw[ultra thick] (1,2)--(2,1);
\draw[ultra thick] (2,1)--(3,0);
\draw[thin] (2,1)--(3,2);
\draw[ultra thick] (3,0)--(4,1);
\draw[ultra thick] (3,2)--(4,1);
\draw[ultra thick] (0,1) arc (180:0:2);

\fill (0,1) circle (2pt);
\fill (1,2) circle (2pt);
\fill (1,0) circle (2pt);
\fill (2,1) circle (2pt);
\fill (2,1) node[above]{$v$};
%\fill (2,-0.5) node[below]{$G$};
\fill (3,0) circle (2pt);
\fill (3,2) circle (2pt);
\fill (4,1) circle (2pt);
\fill (0,0) node[below]{$\phantom{x}$};
\end{tikzpicture}
}

&

\subfigure[$G'$ and $T'$]
{
\begin{tikzpicture}[scale=0.7,style=thick]
\draw[thin] (0,1)--(1,2);
\draw[thin] (0,1)--(1,0);
\draw[ultra thick] (1,0)--(2,0.5);
\draw[ultra thick] (1,2)--(2,1.5);
\draw[ultra thick] (2,0.5)--(3,0);
\draw[thin] (2,1.5)--(3,2);
\draw[ultra thick] (2,1.5)--(3,2);
\draw[ultra thick] (3,0)--(4,1);
\draw[ultra thick] (3,2)--(4,1);
\draw[ultra thick] (0,1) arc (180:0:2);

\fill (0,1) circle (2pt);
\fill (1,2) circle (2pt);
\fill (1,0) circle (2pt);
\fill (2,0.5) circle (2pt);
\fill (2,0.5) node[below]{$v'$};
\fill (2,1.5) circle (2pt);
\fill (2,1.5) node[above]{$v''$};
%\fill (2,-0.5) node[below]{$G'$};
\fill (3,0) circle (2pt);
\fill (3,2) circle (2pt);
\fill (4,1) circle (2pt);
\fill (0,0) node[below]{$\phantom{x}$};
\end{tikzpicture}
}

&

\subfigure[$G''$ and $T''$]
{
\begin{tikzpicture}[scale=0.7,style=thick]
\draw[thin] (0,1)--(1,2);
\draw[ultra thick] (0,1)--(1,2);
\draw[thin] (0,1)--(1,0);
\draw[ultra thick] (1,0)--(2,1);
\draw[ultra thick] (1,2)--(2,1);
\draw[ultra thick] (3,1)--(4,0);
\draw[thin] (3,1)--(4,2);
\draw[ultra thick] (5,1)--(4,0);
\draw[ultra thick] (5,1)--(4,2);
\draw[ultra thick] (0,1) arc (180:0:2.5);

\fill (0,1) circle (2pt);
\fill (1,2) circle (2pt);
\fill (1,0) circle (2pt);
\fill (2,1) circle (2pt);
\fill (2,1) node[above]{$v'$};
%\fill (2.5,-0.5) node[below]{$G'$};
\fill (3,1) circle (2pt);
\fill (3,1) node[above]{$v''$};
\fill (4,0) circle (2pt);
\fill (4,2) circle (2pt);
\fill (5,1) circle (2pt);
\fill (0,0) node[below]{$\phantom{x}$};
\end{tikzpicture}
}
\end{tabular}
\end{center}
\caption{$G'$, $G''$ and $T'$, $T''$ are obtained from $G$ and $T$ using construction from Lemma~\ref{lemma:treeSplittVd}, respectively. Values of their $d$-deficiency are $\xi(G,4,T) = 2$, $\xi(G',4,T') = 1$, and $\xi(G'',4,T'') = 2$. Edges contained in spanning trees are drawn thick.}
\label{fig:example}
\end{figure}

Note that it is not necessary that every graph $G' \in \mathcal{G}$ constructed with splitting procedure as described in Lemma~\ref{lemma:treeSplittVd} has a spanning tree $T'$ such that $\xi(G',d,T') < \xi(G,d,T)$. See Fig.~\ref{fig:example} where construction from Lemma~\ref{lemma:treeSplittVd} is used to obtain graphs $G'$ and $G''$ from a graph $G$ and spanning trees $T'$ and $T''$ from a spanning tree $T$, while $\xi(G,4,T) = 2$, $\xi(G',4,T') = 1$, and $\xi(G'',4,T'') = 2$. 

%%%%%%%%%%%%%%%%%%%%%%%%%%%%%%%%
\section{Proof of Theorem~\ref{thm:main}}
\label{sec:proof}

We now present the proof of Theorem~\ref{thm:main}.

\vspace{10pt}

\proof
Suppose first that the graph $G$ admits an antiparallel $d$-stable trace $W$. Proposition~\ref{prop:dstable} implies that $\delta(G) > d$. If $\Delta(G) < 2d + 2$, then the Fact~\ref{fac:strong-stable} implies that $W$ is also an antiparallel strong trace of $G$. By Theorem~\ref{thm:strong} it then follows that $G$ has a spanning tree $T$ such that each component of the co-tree $G-E(T)$ is even. 

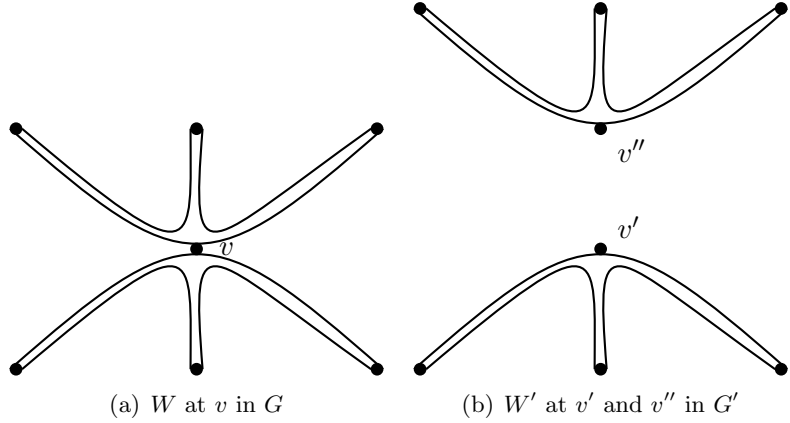
\begin{figure}[ht!]
\begin{center}
\subfigure[$W$ at $v$ in $G$]
{
\begin{tikzpicture}[scale=0.8,style=thick]
\fill (3,2) circle (3pt);
\fill (3.2,2) node[right]{$v$};
\fill (0,0) circle (3pt);
\fill (3,0) circle (3pt);
\fill (6,0) circle (3pt);
\fill (0,4) circle (3pt);
\fill (3,4) circle (3pt);
\fill (6,4) circle (3pt);
\draw (-0.1,4)  .. controls (2.8,1.45) and (3.2,1.45) .. (6.1,4);
\draw (0.1,4)  .. controls (3.1,1.45) and (2.9,2) .. (2.9,4);
\draw (3.1,4)  .. controls (2.9,1.45) and (3.1,2) .. (5.9,4);
\draw (-0.1,0)  .. controls (2.8,2.55) and (3.2,2.55) .. (6.1,0);
\draw (0.1,0)  .. controls (3.1,2.55) and (2.9,2) .. (2.9,0);
\draw (3.1,0)  .. controls (2.9,2.55) and (3.1,2)  .. (5.9,0);
\end{tikzpicture}
}
\subfigure[$W'$ at $v'$ and $v''$ in $G'$]
{
\begin{tikzpicture}[scale=0.8,style=thick]
\fill (3,2) circle (3pt);
\fill (3.1,2) node[above right]{$v'$};
\fill (3,4) circle (3pt);
\fill (3.1,4) node[below right]{$v''$};
\fill (0,0) circle (3pt);
\fill (3,0) circle (3pt);
\fill (6,0) circle (3pt);
\fill (0,6) circle (3pt);
\fill (3,6) circle (3pt);
\fill (6,6) circle (3pt);
\draw (-0.1,6)  .. controls (2.8,3.45) and (3.2,3.45) .. (6.1,6);
\draw (0.1,6)  .. controls (3.1,3.45) and (2.9,4) .. (2.9,6);
\draw (3.1,6)  .. controls (2.9,3.45) and (3.1,4) .. (5.9,6);
\draw (-0.1,0)  .. controls (2.8,2.55) and (3.2,2.55) .. (6.1,0);
\draw (0.1,0)  .. controls (3.1,2.55) and (2.9,2) .. (2.9,0);
\draw (3.1,0)  .. controls (2.9,2.55) and (3.1,2)  .. (5.9,0);
\end{tikzpicture}
}
\end{center}
\caption{Construction from the proof of Theorem~\ref{thm:main}. Note that because $W$ is a double trace there exists a path between the neighbors of $v$ in the upper and the lower part of figure, not shown here. Analogously for $W'$.}
\label{fig:construction3}
\end{figure}

Assume next that $\Delta(G) \geq 2d + 2$. Denote the number of vertices where $W$ has a nontrivial repetitions with $r(W)$. Note that $W$ does not have any nontrivial repetition of order smaller or equal to $d$. We proceed with induction on $r = r(W)$. If $r = 0$, it follows that $W$ is an antiparallel strong trace of $G$ and again, by Theorem~\ref{thm:strong} $G$ has a spanning tree such that each component of the co-tree $G - E(T)$ is even. Let next $r = 1$ and let $v$, $d(v) \geq 2d + 2$, be a unique vertex of $G$ where $W$ has a nontrivial repetition (of order greater than $d$). To make the argument more transparent, assume first that $W$ has exactly two nontrivial repetitions at $v$: $N$ and $N(v) \setminus N$. Obtain graph $G'$ from $G$ with splitting procedure in $v$ using $N$ and $N(v) \setminus N$. Denote new vertices which are adjacent to vertices from $N$ and $N(v) \setminus N$ with $v'$ and $v''$, respectively. See Fig.~\ref{fig:construction3}. Construct double trace $W'$ in $G'$ from $W$ as follows. Let $e = xy$ be an arbitrary (oriented) edge of $W$. If $x \neq v$ and $y \neq v$ we put $e$ into $W'$. If $x = v$ and $y \in N$ then we replace $e$ with $v'y$. Analogously we replace edges where $x = v$ and $y \in N(v) \setminus N$, or $y = v$ and $x \in N$, or $y = v$ and $x \in N(v) \setminus N$ with $v''y$, $xv'$, or $xv''$, respectively. Note first that any edge $e'$ that appears in $G'$ has its unique corresponding edge $e$ in $G$. Since $e$ is traversed twice in the opposite direction in $W$, the edge $e'$ is traversed twice in the opposite direction in $W'$. Hence $W'$ is an antiparallel double trace. For any vertex $u \in V(G) \setminus \{v',v''\}$, $W'$ is without nontrivial repetition at $u$ since otherwise already $W$ would have a nontrivial repetition at $u$. $W'$ is without nontrivial repetitions at $v'$ and $v''$ by construction. Therefore $W'$ is an antiparallel strong trace in $G'$ and Theorem~\ref{thm:santi} the  implies that $G'$ has a spanning tree $T'$ such that each component of $G' - E(T')$ is even. It then follows by Lemma~\ref{lemma:treeContractE} that $G$ has a spanning tree $T$ such that each component of the co-tree $G - E(T)$ is even or contains $v$, which is of degree $\geq 2d + 2$. 

If $W$ has $k > 2$ nontrivial repetitions at $v$, and $N_1, \ldots, N_k \subseteq N(v)$, $\bigcup_{1 \leq i \leq k}N_i = N(v)$ are the sets representing nontrivial repetitions of $W$ at $v$, we construct $G'$ with splitting procedure in $v$ using $N_1, \ldots, N_k$ and denoting new vertices with $v_1, \ldots, v_k$. As above construct double trace $W'$ in $G'$ from $W$, replacing vertex $v$ with new vertex $v_i$ in every edge where $v$ is adjacent to a vertex from $N_i$. Using the same arguments as above it is not difficult to see that $W'$ is an antiparallel strong trace in $G'$ and that $G'$ has a spanning tree $T'$ such that each component of $G' - E(T')$ is even. It then follows by Lemma~\ref{lemma:treeContractE} that $G$ has a spanning tree $T$ such that each component of the co-tree $G - E(T)$ is even or contains $v$, which is of degree $\geq 2d + 2$. We have thus proved that if $G$ has an antiparallel $d$-stable trace with nontrivial repetitions appearing only at vertex $v$ of degree $\geq 2d+2$, then $G$ has a spanning tree $T$ such that each component of the co-tree $G - E(T)$ is even or contains a vertex of degree $\geq 2d + 2$. 

Assume now that $r > 1$ and let $v$ be one of the vertices where $W$ has nontrivial repetitions. We again use the above described construction to obtain graph $G'$ from $G$ by splitting $v$ into $k$ new vertices and finding antiparallel double trace $W'$ in $G'$. It is clear that $r(W') < r(W)$ and hence $G'$ has a spanning tree $T'$ such that each component of the co-tree $G' - E(T')$ is even or contains vertex of degree $\geq 2d + 2$, by induction on $r$. By Lemma~\ref{lemma:treeContractE} it then follows that also $G$ has a spanning tree $T$ such that each component of the co-tree $G - E(T)$ is even or contains a vertex of degree $\geq 2d + 2$ or contains a vertex $v$, which is also of degree $\geq 2d + 2$.

\vspace{10pt}

Conversely, let a graph $G$, $\delta(G) > d$, has a spanning tree $T$ such that each component of the co-tree $G - E(T)$ has an even number of edges or contains a vertex $v$, $d_G(v) \geq 2d + 2$. We proceed by induction on $\xi = \xi(G,2d+2)$. 

If $\xi(G,2d+2) = 0$, then $G$ has a spanning tree $T$ such that each component of the co-tree $G-E(T)$ is even. Theorem~\ref{thm:santi} then implies that $G$ admits an antiparallel strong trace $W$ and since $\delta(G) > d$ it follows from the Fact~\ref{fac:strong-stable} that $W$ is also an antiparallel $d$-stable trace. 

Let $\xi = 1$, $T$ a spanning tree of $G$ which realizes $\xi$, and $v$ a vertex of degree at least $2d + 2$ in a unique odd component of the co-tree $G - E(T)$. Let $\mathcal{G}$ be the family of all graphs obtained from $G$ with splitting procedure in $v$ using two disjoint subsets of $N(v)$ of cardinality $\left\lceil\frac{d(v)}{2}\right\rceil$ and $\left\lfloor\frac{d(v)}{2}\right\rfloor$. Denote new vertices adjacent to $\left\lceil\frac{d(v)}{2}\right\rceil$ and adjacent to $\left\lfloor\frac{d(v)}{2}\right\rfloor$ neighbors of $v$ with $v'$ and $v''$, respectively. By Lemma~\ref{lemma:treeSplittV} there exists a graph $G' \in \mathcal{G}$ such that $G'$ has a spanning tree $T'$ for which each component of the co-tree $G' - E(T')$ is even. Theorem~\ref{thm:santi} then implies that $G'$ admits an antiparallel strong trace $W'$. Construct double trace $W$ in $G$ from $W'$ as follows. Let $e '= xy$ be an arbitrary (oriented) edge of $W'$. If $x,y \notin \{v',v''\}$ we put $e'$ into $W$. If $x \in \{v',v''\}$ then replace $e'$ with $vy$. Analogously we replace edges where $y \in \{v',v''\}$ with $xv$. Note first that any edge $e$ that appears in $G$ has its unique corresponding edge $e'$ in $G'$. Since $e'$ is traversed twice in the opposite direction in $W'$, the edge $e$ is traversed twice in the opposite direction in $W$. Hence $W$ is an antiparallel double trace. For any vertex $u \in V(G) \setminus \{v',v''\}$, $W$ is without nontrivial repetition at $u$ because otherwise already $W'$ would have a nontrivial repetition at $u$. $W$ has exactly two nontrivial repetitions at $v$, $N_{G'}(v')$-repetition and $N_{G'}(v'')$-repetition, by construction. Since $\delta(G) > d$, $|N_{G'}(v')| = \left\lceil\frac{d(v)}{2}\right\rceil > d$, and $|N_{G'}(v')| = \left\lfloor\frac{d(v)}{2}\right\rfloor > d$, $W$ is an antiparallel $d$-stable trace in $G$. We have thus proved that if $G$ has a spanning tree $T$ such that each component of the co-tree $G-E(T)$ is even or contains a vertex of degree at least $2d+2$ and $\xi(G,2d+2) = 1$, then $G$ admits an antiparallel $d$-stable trace. 

Assume now that $\xi > 1$,  $T$ is a spanning tree of $G$ which realizes $\xi$, and $v$is an arbitrary vertex of degree at least $2d + 2$ in an arbitrary odd component of the co-tree $G - E(T)$. As in the case $\xi = 1$, define the family of graphs $\mathcal{G}$. By Lemma~\ref{lemma:treeSplittVd} there exists a graph $G' \in \mathcal{G}$ such that $G'$ has a spanning tree $T'$ for which $\xi(G',2d+2,T') < \xi(G,2d+2,T)$. Since a spanning tree $T$ realizes $(2d+2)$-deficiency of $G$ it follows that $\xi(G',2d+2) < \xi(G,2d+2)$ and hence $G'$ admits an antiparallel $d$-stable trace $W'$ by induction. As above construct an antiparallel double trace $W$ in $G$ from $W'$ by properly replacing $v'$ and $v''$ with $v$. From construction of $W$ it follows that if $W$ has a repetition of order $\leq d$ at some vertex $u \neq v$ already $W'$ has a repetition of order $\leq d$ at vertex $u$, and that if $W$ has a repetition of order $\leq d$ at vertex $v$ also $W'$ has a repetition of order $\leq d$ at $v'$ or $v''$, which is absurd. Therefore $G$ admits an antiparallel $d$-stable trace. 

This concludes the proof of Theorem~\ref{thm:main}.

\qed

\vspace{10pt}

Note that the next corollary of Tutte (1961) and Nash-Williams (1961) tree-packing theorem could also be used to show that in special case every $4$-edge-connected graph $G$, $\delta(G) > d$, with a vertex of degree at least $2d + 2$ or an even Betti number has an antiparallel $d$-stable trace. 

\begin{theorem}
\label{thm:4edgei}
{\rm [Theorem~$2$ of Kundu (1974)]}
Every $4$-edge-connected graph contains two edge-disjoint spanning trees.
\end{theorem}

By Theorem~\ref{thm:4edgei} $4$-edge-connected graph $G$ has a spanning tree $T$ such that a co-tree $G - E(T)$ has exactly one component. 

%%%%%%%%%%%%%%%%%%%%%%%
\section{Concluding remarks}

\subsection{Characterization of graphs admitting double and related traces}

Theorem~\ref{thm:main} also concludes the characterization of graphs admitting (parallel and antiparallel) double, stable, and strong traces, a problem which was first posed by Klav\v zar and Rus (2013). Conditions required for graphs to admit double and related traces are collected in Table~\ref{tab:conditions}. Note again that all the graphs considered are finite and connected. For the sake of clarity of presentation, a spanning tree $T$ of graph $G$ such that each component of the co-tree $G-E(T)$ is even is denoted with $XT$, while a spanning tree $T'$ such that each odd component of the co-tree $G-E(T')$ contains a vertex of degree at least $d$ is denoted with $XT(d)$.

\begin{table}[ht!]
\begin{center}
\begin{tabular}{cc|c|c|c|}
\cline{3-5}
& & \multicolumn{3}{ c| }{$\phantom{X^{X^X}_X}$ CONDITION $\phantom{X^{X^X}_X}$} \\ \cline{3-5}
& & $\phantom{X^{X^X}_X}$ none $\phantom{X^{X^X}_X}$ & parallel & antiparallel \\ \cline{1-5}
\multicolumn{1}{ |c| }{\multirow{3}{*}{\rotatebox{90}{\mbox{$\phantom{X^{X^X}}$ TRACE $\phantom{X}$}}}} &
\multicolumn{1}{ |c| }{double} & 
$\displaystyle{{\phantom{X^{X^X}}\forall G \phantom{X^{X^X}}\atop\text{{\rm [K\" onig (1936)]}}}}$ & 
$\displaystyle{{\phantom{X}{\rm Eulerian} \phantom{X}\atop\text{{\rm [Theorem~2 of Vestergaard (1975)]}}}}$ & 
$\displaystyle{{\forall G \atop \text{{\rm [Tarry (1895)]}}}}$ 
\\ \cline{2-5}
\multicolumn{1}{ |c  }{}                        &
\multicolumn{1}{ |c| }{$d$-stable} & 
$\displaystyle{{\phantom{X^{X^X}}\delta(G) > d \phantom{X^{X^X}}\atop \text{{\rm [Theorem~3.4 of (*)]}}}}$ & 
$\displaystyle{{\phantom{X}\delta(G) > d \wedge {\rm Eulerian} \phantom{X}\atop \text{{\rm [Theorem~5.4 of (*)]}}}}$ & 
$\displaystyle{{\delta(G) > d \wedge \exists {\rm XT}(2d+2) \atop {\rm (Theorem~\ref{thm:main})}}}$ 
\\ \cline{2-5}
\multicolumn{1}{ |c  }{}                        &
\multicolumn{1}{ |c| }{strong} & 
$\displaystyle{{\phantom{X^{X^X}}\forall G \phantom{X^{X^X}}\atop \text{{\rm [Theorem~2.4 of (*)]}}}}$ & 
$\displaystyle{{\phantom{X}{\rm Eulerian} \phantom{X}\atop \text{{\rm [Theorem~5.3 of (*)]}}}}$ & 
$\displaystyle{{\exists {\rm XT} \atop \text{{\rm [Theorem~4.1 (*)]}}}}$ 
\\ \cline{1-5}
\end{tabular}
\end{center}
\caption{Required conditions for graphs to admit presented double traces. Note that Theorems denoted with asterisk were proven by Fijav\v z et al. (2014).}
\label{tab:conditions}
\end{table}

\subsection{Polynomial algorithm for $d = 1$}

Furst et al. (1988) presented a polynomial time algorithm which for an arbitrary graph $G$ finds a spanning tree $T$ such that that each component of the co-tree $G-E(T)$ is even or determine that no such $T$ exists. It uses matroids parity (or matching). Thomassen (1990) modified the algorithm so it also works for a spanning tree $T$ such that each component of the co-tree $G-E(T)$ is even or contains a vertex of degree at least $4$. It seems that the algorithm by Furst et al. (1988) modified in a similar way as by Thomassen (1990) would also work for spanning trees, which satisfy conclusion from Theorem~\ref{thm:main} for arbitrary $d$.  

\subsection{Pseudosurfaces}

We present an alternative characterization which may help to visualize the graphs considered in Theorem~\ref{thm:main}. Fijav\v z et al. (2014) observed that graph $G$ admits a strong trace if and only if $G$ has a single face embedding in some surface, and that graph $G$ admits an antiparallel strong trace if and only if $G$ has a single face embedding in some orientable surface. Characterizations of graphs admitting single face embeddings presented by Edmonds (1965) and Behzad et al. (1979) were then used to characterize graphs which admit (antiparallel) strong traces. Here we can do the reverse and use the characterization of graphs admitting antiparallel $d$-stable trace to characterize graphs which have some other properties as well.

A {\em pinched open disk} is defined by Pisanski and Poto\v cnik (2003) as a topological space obtained from $k$ copies of open disks by identifying their respective centers to a single vertex. {\em Pseudosurfaces} are then obtained with relaxation of the definition of surfaces by allowing the neighborhoods of points to be homeomorphic not only to open disks or half-disks but also to pinched open disks. A pseudosurface can be obtained from a single surface by identifying its points, from two or more surfaces by pairwise identifying some number of points on one sphere with points on the other surfaces, or by combination of previous two operations. A pseudosurface is an {\em orientable pseudosurface} if and only if all the surfaces used in its construction are orientable. Since single face embeddings can only be defined for the first of the three presented constructions of pseudosurfaces, we pose next observation only for those pseudosurfaces. Theorem~\ref{thm:main} actually characterizes graphs which have a single face embedding in some orientable pseudosurface obtained from a single surface by identifying its points, such that for every vertex $v$ and every neighborhood $N$ of $v$ homeomorphic to an open disk, $v$ has at least $d+1$ neighbors in $N$. Note that if $v$ is a vertex where $k$ copies of open disks were identified and form a pinched open disks, then $v$ has $k$ disjoint neighborhoods homeomorphic to an open disk.  

This alternative characterization is especially interesting since embeddings of graphs in pseudosurfaces have received much less attention than embeddings in surfaces, although the literature does contain some results (Petroelje 1971; Abrams and Slilaty~2006).

\subsection{Bionanostructures}

\begin{figure}[ht!]
\begin{center}
\begin{tikzpicture}[scale=1.0,style=thick]
\fill (3,2) circle (3pt);
\fill (3.2,2) node[right]{$v$};
\fill (0,0) circle (3pt);
\fill (3,0) circle (3pt);
\fill (6,0) circle (3pt);
\fill (0,4) circle (3pt);
\fill (3,4) circle (3pt);
\fill (6,4) circle (3pt);
\draw (-0.1,4)  .. controls (2.8,1.45) and (3.2,1.45) .. (6.1,4);
\draw (0.1,4)  .. controls (3.1,1.45) and (2.9,2) .. (2.9,4);
\draw (3.1,4)  .. controls (2.9,1.45) and (3.1,2) .. (5.9,4);
\draw (-0.1,0)  .. controls (2.8,2.55) and (3.2,2.55) .. (6.1,0);
\draw (0.1,0)  .. controls (3.1,2.55) and (2.9,2) .. (2.9,0);
\draw (3.1,0)  .. controls (2.9,2.55) and (3.1,2)  .. (5.9,0);
\end{tikzpicture}
\end{center}
\caption{A $3$-repetition in a vertex $v$ of degree $6$}
\label{fig:repetition}
\end{figure}
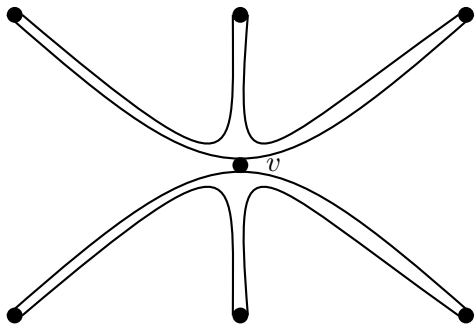

It was mentioned in the Introduction that every polyhedron $P$ which is composed from a single polymer chain can be naturally represented with a graph $G(P)$ of the polyhedron. Since every edge of such a graph $G(P)$ corresponds to exactly two segments of polymer chain, single polymer chain corresponds to a double trace in $G(P)$.
Furthermore, parallel edges (of double traces) represent pairs of two coiled-coil-forming segments aligned in the same direction while antiparallel edges represent pairs of two coiled-coil-forming segments aligned in the opposite direction. With the recent demonstration by Ko\v car et al. (2016) of DNA-based polyhedra assembled from a single chain and the fact that DNA segments are always aligned in opposite direction, antiparallel double traces became the focus of our research.
At last, a polyhedron $P$ composed from a single polymer chain corresponding to a stable trace with repetitions may fold to a polyhedron different from $P$, as a vertex may indeed split into a collection of independent vertices of smaller degree, see also Fig.~\ref{fig:repetition}. Therefore, $d$-stable traces of higher order (higher $d$) provide more stable structures.

In conclusion the theoretical framework established in this contribution will contribute to the design of new bionanostructures composed on natural building blocks.

\section*{Acknowledgements}

The author is grateful to Roman Jerala for help with biological aspect of the problem. Many helpful suggestions were also given by Dan Archdeacon, Ga\v sper Fijav\v z, Luis Goddyn, and Sandi Klav\v zar.

%%%%%%%%%%%%%%%%%%%%%%%%%%%

\end{document}